\documentclass[11pt]{article}
\usepackage{amsmath}
\usepackage{amssymb}
\usepackage{natbib}

\begin{document}
\bibliographystyle{plainnat}
\setcitestyle{numbers,square}

\title{EXPLICIT INVERSION FOR TWO \\
       BROWNIAN--TYPE MATRICES}

\author{Florendia Valvi$^1$, Vassilis Geroyannis$^2$ \\
        $^1$Department of Mathematics, University of Patras, Greece \\
        $^2$Department of Physics, University of Patras, Greece \\
        Email: $^1$fvalvi@upatras.gr, $^2$vgeroyan@upatras.gr}
\maketitle

\begin{abstract}
We present explicit inverses of two Brownian--type matrices, which are defined as Hadamard products of certain already known matrices. The matrices under consideration are defined by $3n\!\!-\!\!1$ parameters and their lower Hessenberg form inverses are expressed analytically in terms of these parameters. Such matrices are useful in the theory of digital signal processing and in testing matrix inversion algorithms. \\   
\\
\textbf{Keywords:} Brownian Matrix; Hadamard Product; Hessenberg Matrix; Numerical Complexity; Test Matrix
\end{abstract}

\section{Introduction}
Brownian matrices are frequently involved in problems concerning ``digital signal processing''. In particular, Brownian motion is one of the most common linear models used for representing nonstationary signals. The covariance matrix of a discrete--time Brownian motion has, in turn, a very characteristic structure, the so-called ``Brownian matrix''.  

In \citep{her69} (Eq.~(2)) the explicit inverse of a class of matrices $G_n = [\beta_{ij}]$ with elements
\begin{equation}
\beta_{ij} = \left\{ \begin{array}{ll} b_j, & \ \ i \leqslant j, \\
                                       a_j, & \ \ i > j.
                     \end{array}
             \right.
\end{equation}
is given. On the other hand, the analytic expressions of the inverses of two symmetric matrices $K = [\kappa_{ij}]$ and $N = [\nu_{ij}]$, where
\begin{equation}
\kappa_{ij} = k_i \quad \textrm{and} \quad \nu_{ij} = k_j, \ \ i \leqslant j,
\end{equation}  
respectively, are presented in \citep{val77} (first equation in p.~113, and Eq.~(1), respectively). The matrix $K$ is a special case of Brownian matrix and $G_n$ is a lower Brownian matrix, as they have been defined in \citep{gba86} (Eq.~(2.1)). Earlier, in \citep{pic83} (paragraph following Eq.~(3.3))  the term ``pure Brownian matrix'' for the type of the matrix $K$ has introduced. Furthermore, in \citep{ckm82} (discussion concerning Eqs.~(28)--(30)) the so-called ``diagonal innovation matrices'' (DIM) have been treated, special cases of which are the matrices $K$ and $N$.  

In the present paper, we consider two matrices $A_1$ and $A_2$ defined by 
\begin{equation}
A_1 = K \circ G_n \ \ \ \textrm{and} \ \ \ A_2 = N \circ G_n,
\label{HP}
\end{equation}
where the symbol $\circ$ denotes the Hadamard product. Hence, the matrices have the forms
\begin{equation}
A_1 = \left[
      \begin{array}{cccccc}
             k_1 b_1 & k_1 b_2 & k_1 b_3 & \ldots & k_1 b_{n-1}     & k_1 b_n     \\
             k_1 a_1 & k_2 b_2 & k_2 b_3 & \ldots & k_2 b_{n-1}     & k_2 b_n     \\
             k_1 a_1 & k_2 a_2 & k_3 b_3 & \ldots & k_3 b_{n-1}     & k_3 b_n     \\
             \ldots  &         &         &        &                 &             \\
             k_1 a_1 & k_2 a_2 & k_3 a_3 & \ldots & k_{n-1} b_{n-1} & k_{n-1} b_n \\
             k_1 a_1 & k_2 a_2 & k_3 a_3 & \ldots & k_{n-1} a_{n-1} & k_n b_n
      \end{array}
      \right]
\label{thematrixA1}
\end{equation} 
and
\begin{equation}
A_2 = \left[
      \begin{array}{cccccc}
             k_1 b_1 & k_2 b_2 & k_3 b_3 & \ldots & k_{n-1} b_{n-1}     & k_n b_n         \\
             k_2 a_1 & k_2 b_2 & k_3 b_3 & \ldots & k_{n-1} b_{n-1}     & k_n b_n         \\
             k_3 a_1 & k_3 a_2 & k_3 b_3 & \ldots & k_{n-1} b_{n-1}     & k_n b_n         \\
             \ldots  &         &         &        &                     &                 \\
             k_{n-1} a_1 & k_{n-1} a_2 & k_{n-1} a_3 & \ldots & k_{n-1} b_{n-1} & k_n b_n \\
             k_n a_1     & k_n a_2     & k_n a_3     & \ldots & k_n a_{n-1}     & k_n b_n
      \end{array}
      \right].
\label{thematrixA2}
\end{equation}

Let us now define for a matrix $B=[b_{ij}]$ the terms ``pure upper Brownian matrix'' and ``pure lower Brownian matrix'', for the elements of which the following relations are respectively valid
\begin{equation}
b_{i,j+1}=b_{ij}, \quad i \leqslant j, \qquad \mathrm{and} \qquad 
                  b_{i+1,j}=b_{ij}, \quad i \geqslant j.
\end{equation}

The matrix $A_1$ (Eq.~(\ref{thematrixA1})) is a lower Brownian matrix. Furthermore, the matrix $P \, N \, P$, where $P=[p_{ij}]$ is the permutation matrix with elements 
\begin{equation}
p_{ij} = \left\{ \begin{array}{ll} 1, & \ \ i+j=n+1, \\
                                       0, & \ \ \mathrm{otherwise},
                     \end{array}
             \right.
\end{equation}
is a pure Brownian matrix and $P \, G_n \, P$ a pure lower Brownian matrix. Hence, their Hadamard product $(P \, N \, P) \circ (P \, G_n \, P)$ gives a pure lower Brownian matrix, that is, the matrix $P \, A_2 \, P$. 

In the following sections, we deduce in analytic form the inverses and determinants of the matrices $A_1$ and $A_2$; and we study the numerical complexity on evaluating $A_1^{-1}$ and $A_2^{-1}$.

\section{The Inverse and Determinant of $A_1$}
\label{A1proof}
The inverse of $A_1$ is a lower Hessenberg matrix expressed analytically by the $3n\!\!-\!\!1$ parameters defining $A_1$. In particular, the inverse $A_1^{-1} = [\alpha_{ij}]$ has elements given by the relations
\begin{equation}
\alpha_{ij} = \left\{ 
              \begin{array}{ll} 
              \frac{\displaystyle k_{i+1} b_{i-1} - k_{i-1} a_{i-1}}{\displaystyle c_{i-1} c_i}, & \ \ \ i = j \neq 1,n, \\
                                                                     &                       \\
              \frac{\displaystyle k_2}{\displaystyle k_1 c_1}, & \ \ \ i = j = 1,                                        \\
                                                                     &                       \\
              \frac{\displaystyle b_{n-1}}{\displaystyle c_{n-1} c_n}, & \ \ \ i = j = n,                                \\
                                                                     &                       \\
              (-1)^{i+j} \frac{\displaystyle d_{j-1} g_i \prod_{\nu=j+1}^{i-1} k_{\nu} f_{\nu}}
                              {\displaystyle \prod_{\nu=j-1}^i c_{\nu}} \ , 
                                                           & \ \ \ i - j \geqslant 1,        \\
                                                                     &                       \\
              - \ \frac{\displaystyle 1}{\displaystyle c_i}, & \ \ \ j - i = 1,                                          \\
                                                                     &                       \\
              0, & \ \ \ j - i > 1,
              \end{array}
              \right.
\label{A1case}    
\end{equation}
where 
\begin{equation}
\left\{
\begin{array}{llll}
c_i = k_{i+1} b_i - k_i a_i,                 & \ \ i = 1,2,\ldots,n-1, & \ \ c_0 =1,     & \ \ c_n = b_n,
\\
d_i = k_{i+1} a_{i+1} b_i  - k_i a_i b_{i+1}, & \ \ i = 1,2,\ldots,n-2, & \ \ d_0 = a_1, &
\\
f_i = a_i - b_i,                             & \ \ i=2,3,\ldots,n-1,   &                 &
\\
g_i = k_{i+1} - k_i,                         & \ \ i=2,3,\ldots,n-1,   & \ \ g_n = 1,    &
\end{array}
\right.
\label{A1conventions}
\end{equation}
with 
\begin{equation}
\prod_{\nu=j+1}^{i-1} k_{\nu} f_{\nu} = 1 \ \ \ \textrm{if} \ \ \ i = j + 1,
\label{Product}
\end{equation}
and with the obvious assumptions
\begin{equation}
k_1 \neq 0 \ \ \ \ \textrm{and} \ \ \ \ c_i \neq 0, \ \ \ i=1,2,\ldots,n.
\end{equation}

To prove that the relations~(\ref{A1case})--(\ref{Product}) give the inverse matrix $A_1^{-1}$, we  reduce $A_1$ to the identity matrix $I$ by applying a number of elementary row transformations. Then the product of the corresponding elementary matrices gives the inverse matrix of $A_1$. These transformations are defined by the following sequence of row operations. \\ 
\leftline{ }
\textbf{Operation~1} (applied on $A_1$ and on the identity matrix $I$): 
\begin{equation*}
\textrm{row} \ i - \frac{k_i}{k_{i+1}} \times \textrm{row} \ (i+1), \ \ i=1,2,\ldots,n-1,
\end{equation*}
which transforms $A_1$ into the lower triangular matrix $C_1$ given by
\begin{equation*}
\left[ \begin{array}{llllll}
       \frac{k_1 \left( k_2 b_1 - k_1 a_1 \right)}{k_2} & 0 & 0 & \ldots & 0 & 0 \\
       & & & & &                                                                 \\
       \frac{k_1 a_1 \left( k_3 - k_2 \right)}{k_3} & \frac{k_2 \left( k_3 b_2 - k_2 a_2 \right)}{k_3}
                                                            & 0 & \ldots & 0 & 0 \\
       & & & & &                                                                 \\
       \frac{k_1 a_1 \left( k_4 -k_3 \right)}{k_4} & \frac{k_2 a_2 \left( k_4 -k_3 \right)}{k_4} 
             & \frac{k_3 \left( k_4 b_3 - k_3 a_3 \right)}{k_4} & \ldots & 0 & 0 \\
       \ldots & & & & &                                                          \\
       \frac{k_1 a_1 \left( k_n - k_{n-1} \right)}{k_n} & \frac{k_2 a_2 \left( k_n - k_{n-1} \right)}{k_n} &
       \frac{k_3 a_3 \left( k_n - k_{n-1} \right)}{k_n} & \ldots & \frac{k_{n-1} \left( k_n b_{n-1} - k_{n-1} a_{n-1}
       \right)}{k_n} & 0                                                               \\
       & & & & &                                                                 \\
       k_1 a_1 & k_2 a_2 & k_3 a_3 & \ldots & k_{n-1} a_{n-1} & k_n b_n
       \end{array}
       \right],
\end{equation*} 
and the identity matrix $I$ into the upper bidiagonal matrix $F_1$ with main diagonal
\begin{displaymath}
\left( 1 \ , \ 1 \ , \ \ldots \ , \ 1 \right)
\end{displaymath}
and upper first diagonal
\begin{displaymath}
\left( -\frac{k_1}{k_2} \ , \ -\frac{k_2}{k_3} \ , \ \ldots \ , \ -\frac{k_{n-1}}{k_n} \right).
\end{displaymath} 
\leftline{ }
\textbf{Operation~2} (applied on $C_1$ and $F_1$): 
\begin{equation*}
\textrm{row} \ i - \ \frac{k_i g_i}{k_{i+1} g_{i-1}} \ \times
              \textrm{row} \ (i-1), \ \ i=n,n-1,\ldots,3, \ \ \ k_{n+1} = 1,
\end{equation*}
which derives a lower bidiagonal matrix $C_2$ with main diagonal 
\begin{displaymath}
\left( \frac{k_1 c_1}{k_2} \ , \ \frac{k_2 c_2}{k_3} \ , \ \ldots \ , \ \frac{k_{n-1} c_{n-1}}{k_n} \ , \ k_n c_n
\right)
\end{displaymath}
and lower first diagonal
\begin{displaymath}
\left( 
\frac{k_1 a_1 g_2}{k_3} \ , \ \frac{k_2 k_3 g_3 f_2}{k_4 g_2} \ , \ \ldots \ , \ \frac{k_{n-2} k_{n-1} g_{n-1} f_{n-2}}{k_n g_{n-2}} \ , \ \frac{k_{n-1} k_n f_{n-1}}{g_{n-1}}
\right);
\end{displaymath} 
while the matrix $F_1$ is transformed into the tridiagonal matrix $F_2$ given by  
\begin{displaymath}
\left[
\begin{array}{cccccc} 
       1 & -\frac{k_1}{k_2} & 0 & \ldots & 0 & 0 \\
       0 & 1 & -\frac{k_2}{k_3} & \ldots & 0 & 0 \\
       \ \ 0 \ \ & \ \ -\frac{k_3 g_3}{k_4 g_2} \ \ & \ \ 1+\frac{k_2 g_3}{k_4 g_2} \ \ & \ \ \ldots \ \ & \ \ 0 \ \ & \ \ 0 \ \ \\
       \ldots & & & & &  \\
       \ \ 0 \ \ & \ \ 0 \ \ & \ \ 0 \ \ & \ \ \ldots \ \ & \ \ 1+\frac{k_{n-2} g_{n-1}}{k_n g_{n-2}} \ \ & \ \ -\frac{k_{n-1}}{k_n} \ \ \\
       0 & 0 & 0 & \ldots & -\frac{k_n}{g_{n-1}} & 1+\frac{k_{n-1}}{g_{n-1}} 
\end{array}
\right].
\end{displaymath}
\leftline{ }
\textbf{Operation~3} (applied on $C_2$ and $F_2$): 
\begin{equation*}
\textrm{row} \ 2 - \ \frac{k_2 a_1 g_2}{k_3 c_1} \ \times \textrm{row} \ 1 \quad \textrm{and}
    \quad 
    \textrm{row} \ i -  
    k_i \ \frac{k_i g_i f_{i-1}}{k_{i+1} g_{i-1} c_{i-1}} \ \times 
    \textrm{row} \ (i-1), \ \ i=3,4,\ldots,n,
\end{equation*} 
which derives the diagonal matrix 
\begin{displaymath}
C_3=
\left\lceil \frac{k_1 c_1}{k_2} \ \ \ \frac{k_2 c_2}{k_3} \ \ \ \ldots \ \ \ \frac{k_{n-1} c_{n-1}}{k_n} \ \ \ k_n c_n \right\rfloor,
\end{displaymath} 
and, respectively, the lower Hessenberg matrix $F_3$ given by  
\begin{equation*}
\left[
\begin{array}{ccccc}
       1 & -\frac{k_1}{k_2} & \ldots & 0 & 0                                                              \\
       & & & & 
\\
       -\frac{k_2 a_1 g_2}{k_3 c_0 c_1} & \frac{k_2 \left( k_3 b_1 - k_1 a_1 \right)}{k_3 c_1} & 
                                                                                \ldots & 0 & 0            \\
       & & & & 
\\
       \frac{k_3 a_1 g_3 k_2 f_2}{k_4 c_0 c_1 c_2} & -\frac{k_3 d_1 g_3}{k_4 c_1 c_2} & \ldots & 0 & 0    \\
       \ldots & & & &                                                                                     \\
       \frac{s \, k_{n-1} a_1 g_{n-1} k_2 f_2 \ldots k_{n-2} f_{n-2}}{k_n c_0 c_1 \ldots c_{n-2}} &
       \frac{s \, k_{n-1} d_1 g_{n-1} k_3 f_3 \ldots k_{n-2} f_{n-2}}{k_n c_1 c_2 \ldots c_{n-2}} & \ldots &
       \frac{k_{n-1} \left( k_n b_{n-2} - k_{n-2} a_{n-2} \right)}{k_n c_{n-2}} & - \frac{k_{n-1}}{k_n}   \\
       & & & & 
\\
       \frac{s \, k_n a_1 g_n k_2 f_2 \ldots k_{n-1} f_{n-1}}{c_0 c_1 \ldots c_{n-1}} & 
       \frac{s \, k_n d_1 g_n k_3 f_3 \ldots k_{n-1} f_{n-1}}{c_1 c_2 \ldots c_{n-1}} & \ldots &
       -\frac{k_n d_{n-2}}{c_{n-2} c_{n-1}} & \frac{k_n b_{n-1}}{c_{n-1}}
\end{array}
\right],
\end{equation*}
with the symbol $s$ standing for the quantity $(-1)^{i+j}$. \\
\leftline{ } 
\textbf{Operation~4} (applied on $C_3$ and $F_3$): 
\begin{equation*}
\frac{k_{i+1}}{k_i c_i} \times \textrm{row} \ i, \ \ i=1,2,\ldots,n,
\end{equation*}
which transforms $C_3$ into the identity matrix $I$ and the matrix $F_3$ into the inverse $A_1^{-1}$. 

The determinant of $A_1$ takes the form 
\begin{equation}
\textrm{det}\left(A_1\right) = k_1 b_n \left( k_2 b_1 - k_1 a_1 \right)
            \left( k_3 b_2 - k_2 a_2 \right) \ \ldots \ 
            \left( k_n b_{n-1} - k_{n-1} a_{n-1} \right).
\end{equation}
Evidently, $A_1$ is singular if $k_1 = 0$ or, considering the relations~(\ref{A1conventions}), if $c_i = 0$ for some $i \in \left\{ 1,2,\ldots,n \right\}$.

\section{The Inverse and Determinant of $A_2$}
In the case of $A_2$, its inverse $A_2^{-1} = [\alpha_{ij}]$ is a lower Hessenberg matrix with elements given by the relations 
\begin{equation}
\alpha_{ij} = \left\{ 
              \begin{array}{ll} 
              \frac{\displaystyle k_{i-1} b_{i-1} - k_{i+1} a_{i-1}}
                   {\displaystyle c_{i-1} c_i}, & \ \ \ i = j \neq 1,n,                      \\
                                                                     &                       \\
              \frac{\displaystyle 1}{\displaystyle c_1}, & \ \ \ i = j = 1,                  \\
                                                                     &                       \\
              \frac{\displaystyle k_{n-1} b_{n-1}}
                   {\displaystyle k_n c_{n-1} c_n}, & \ \ \                i = j = n,        \\
                                                                     &                       \\
              (-1)^{i+j} \frac{\displaystyle d_{j-1} g_i \prod_{\nu=j+1}^{i-1} k_{\nu} f_{\nu}}
                              {\displaystyle \prod_{\nu=j-1}^i c_{\nu}}, 
                                                           & \ \ \ i - j \geqslant 1,        \\
                                                                     &                       \\
              - \ \frac{\displaystyle 1}{\displaystyle c_i}, & \ \ \ j - i = 1,              \\
                                                                     &                       \\
              0, & \ \ \ j - i > 1,
              \end{array}
              \right. 
\label{A2case}   
\end{equation}
where 
\begin{equation}
\left\{
\begin{array}{llll}
c_i = k_i b_i - k_{i+1} a_i,                 & \ \ i = 1,2,\ldots,n-1, & \ \ c_0 =1,    & \ \ c_n = b_n,
\\
d_i = k_i a_{i+1} b_i - k_{i+1} a_i b_{i+1}, & \ \ i = 1,2,\ldots,n-2, & \ \ d_0 = a_1, &
\\
f_i = a_i - b_i,                             & \ \ i=2,3,\ldots,n-1,   &                &
\\
g_i = k_i - k_{i+1},                         & \ \ i=2,3,\ldots,n-1,   & \ \ g_n = 1,   & 
\end{array}
\right.
\label{A2conventions}
\end{equation}
with
\begin{equation}
\prod_{\nu=j+1}^{i-1} k_{\nu} f_{\nu} = 1 \ \ \ \textrm{if} \ \ \ i = j + 1,
\label{Product2}
\end{equation}
and with the obvious assumptions
\begin{equation}
k_n \neq 0 \ \ \ \ \textrm{and} \ \ \ \ c_i \neq 0, \ \ \ i=1,2,\ldots,n.
\end{equation}

In order to prove that the relations~(\ref{A2case})--(\ref{Product2}) give the inverse matrix $A_2^{-1}$, we follow a similar manner to that of Sec.~\ref{A1proof}. \\
\leftline{ }
\textbf{Operation~1} (applied on $A_2$ and on the identity matrix $I$):
\begin{equation*}
\textrm{row} \ i - \textrm{row} \ (i+1), \ \ \ i = 1,2,\ldots,n-1,
\end{equation*}
which transforms $A_2$ into the lower triangular matrix $D_1$ equal to 
\begin{equation*}
\left[
\begin{array}{ccccc}
       k_1 b_1 - k_2 a_1 & 0 & \ldots & 0 & 0 \\
       a_1 \left( k_2 - k_3 \right) & k_2 b_2 - k_3 a_2 & \ldots & 0 & 0 \\
       \ldots & & & &  \\
       \ \ a_1 \left( k_{n-1} - k_n \right) \ \ & \ \ a_2 \left( k_{n-1} - k_n \right) \ \ & \ \ \ldots \ \ & 
       \ \ k_{n-1} b_{n-1} - k_n a_{n-1} \ \ & \ \ 0 \ \ \\
       k_n a_1 & k_n a_2 & \ldots & k_n a_{n-1} & k_n b_n
\end{array}
\right],
\end{equation*} 
and the identity matrix $I$ into the bidiagonal matrix $L_1$ 
with main diagonal
\begin{displaymath}
\left( 1 \ , \ 1 \ , \ \ldots \ , \ 1 \ , \ 1 \right)
\end{displaymath}
and upper first diagonal
\begin{displaymath}
\left( -1 \ , \ -1 \ , \ \ldots \ , \ -1 \ , \ -1 \right).
\end{displaymath} 
\leftline{ }
\textbf{Operation~2} (applied on $D_1$ and $L_1$):
\begin{equation*}
\textrm{row} \ n - \frac{k_n}{g_{n-1}} \times \textrm{row} \ (n-1) \ \ \ \textrm{and} \ \ \ 
\textrm{row} \ i - \frac{g_i}{g_{i-1}} \times \textrm{row} \ (i-1), \ \ \ i=n-1,n-2,\ldots,3,
\end{equation*}
which derives the lower bidiagonal matrix $D_2$ with main diagonal
\begin{displaymath}
\left( c_1 \ , \ c_2 \ , \ \ldots \ , \ c_{n-1} \ , \ k_n c_n \right)
\end{displaymath}
and lower first diagonal
\begin{displaymath}
\left( a_1 g_2 \ , \ \frac{g_3 k_2 f_2}{g_2} \ , \ \ldots \ , \ \frac{g_{n-1} k_{n-2} f_{n-2}}
                          {g_{n-2}} \ , \ \frac{k_n k_{n-1} f_{n-1}}{g_{n-1}} \right);
\end{displaymath} 
while the matrix $L_1$ is transformed into the tridiagonal matrix $L_2$ with main diagonal
\begin{displaymath}
\left( 1 \ , \ 1 \ , \ 1 + \frac{g_3}{g_2} \ , \ \ldots \ , \ 1 + \frac{g_{n-1}}{g_{n-2}} \ , \
       1 + \frac{k_n}{g_{n-1}} \right),
\end{displaymath}
upper first diagonal
\begin{displaymath}
\left( -1 \ , \ -1 \ , \ \ldots \ , \ -1 \ , \ -1 \right)
\end{displaymath}
and lower first diagonal
\begin{displaymath}
\left( 0 \ , \ -\frac{g_3}{g_2} \ , \ \ldots \ , \ -\frac{g_{n-1}}{g_{n-2}} \ , \
       -\frac{k_n}{g_{n-1}} \right).
\end{displaymath}
\leftline{ }
\textbf{Operation 3} (applied on $D_2$ and $L_2$):
\begin{equation*}
\textrm{row} 2 - \frac{a_1 g_2}{c_1}\textrm{row} 1,  \ \ \textrm{row} i - \frac{g_i k_{i-1} f_{i-1}}{g_{i-1} c_{i-1}}\textrm{row} (i\!-\!1), \ \ \ldots, \ \
        \textrm{row} n - \frac{k_n k_{n-1} f_{n-1}}{g_{n-1} c_{n-1}}\textrm{row} (n\!-\!1),
\end{equation*}
with $i=3,4,\ldots,n\!-\!1$, which yields the diagonal matrix $D_3$, 
\begin{displaymath}
D_3=\left\lceil c_1 \ \ \ c_2 \ \ \ \ldots \ \ \ c_{n-1} \ \ \ k_n c_n \right\rfloor,
\end{displaymath} 
and the lower Hessenberg matrix $L_3$ equal to
\begin{displaymath}
\left[
\begin{array}{ccccc}
       1 & -1 & \ldots & 0 & 0                                                                 \\
       & & & & \\                 
       -\frac{a_1 g_2}{c_0 c_1} & \frac{k_1 b_1 - k_3 a_1}{c_1} & \ldots & 0 & 0               \\
       & & & & \\ 
       \frac{a_1 g_3 k_2 f_2}{c_0 c_1 c_2} & -\frac{d_1 g_3}{c_1 c_2} & \ldots & 0 & 0         \\
       \ldots & & & &                                                                          \\
       \frac{s \, a_1 g_{n-1} k_2 f_2 \ldots k_{n-2} f_{n-2}}{c_0 c_1 \ldots c_{n-2}} &
       \frac{s \, d_1 g_{n-1} k_3 f_3 \ldots k_{n-2} f_{n-2}}{c_1 \ldots c_{n-2}}     &
       \ldots                                                                         &
       \frac{k_{n-2} b_{n-2} - k_n a_{n-2}}{c_{n-2}}                                  & 
       -1                                                                                       \\
       & & & & \\ 
       \frac{s \, k_n a_1 g_n k_2 f_2 \ldots k_{n-1} f_{n-1}}{c_0 c_1 c_2 \ldots c_{n-1}}    &
       \frac{s \, k_n d_1 g_n k_3 f_3 \ldots k_{n-1} f_{n-1}}{c_1 c_2 \ldots c_{n-1}}        &
       \ldots                                                                                &
      -\frac{k_n d_{n-2} g_n}{c_{n-2} c_{n-1}} & \frac{k_{n-1} b_{n-1}}{c_{n-1}}          
\end{array}
\right],
\end{displaymath}
where the symbol $s$ stands for $(-1)^{i+j}$. \\
\leftline{ }
\textbf{Operation 4} (applied on $D_3$ and $L_3$):
\begin{equation*}
\frac{1}{c_i} \times \textrm{row} \ i, \ \ i=1,2,\ldots,n-1, \ \ \ \ \textrm{and} \ \ \
   \frac{1}{k_n c_n} \times \textrm{row} \ n,
\end{equation*} 
which transforms $D_3$ into the identity matrix $I$ and $L_3$ into the inverse $A_2^{-1}$.

The determinant of $A_2$ has the form 
\begin{equation}
\textrm{det}(A_2) = k_n b_n \left( k_1 b_1 - k_2 a_1 \right) \left( k_2 b_2 - k_3 a_2 \right)
                    \ldots \left( k_{n-1} b_{n-1} - k_n a_{n-1} \right),
\end{equation}
which shows in turn that the matrix $A_2$ is singular if $k_n = 0$, or, adopting the conventions (\ref{A2conventions}), if $c_i = 0$ for some $i \in \left\{ 1,2,\ldots,n \right\}$.

\section{Numerical Complexity}
\label{complexity}
The relations~(\ref{A1case}) and (\ref{A2case}) lead to recurrence formulae, by which the inverses $A_1^{-1}$ and $A_2^{-1}$, respectively, are computed in $O(n^2)$ multiplications/divisions and $O(n)$ additions/substractions. In fact, the recursive algorithm
\begin{equation}
\alpha_{i,i+1} = - 1 / c_i, \ \ \ i=1,2,\ldots,n-1,
\label{fra1}
\end{equation}
\begin{equation}
\alpha_{ii} = - \alpha_{i,i+1} + \frac{b_{i-1} g_i}{c_{i-1} c_i}, \ \ \ 
              i=2,3,\ldots,n-1, \ \ \ \alpha_{11}=\frac{k_2}{k_1 c_1}, \ \ \ 
              \alpha_{nn} = \frac{b_{n-1}}{c_{n-1} c_n},
\end{equation}
\begin{equation}
\alpha_{i,i-1} = - \frac{d_{i-2} g_i}{c_{i-2} c_{i-1} c_i}, \ \ \ i=2,3,\ldots,n,
\end{equation}
\begin{equation}
\alpha_{i,i-s-1} = - \ \frac{d_{i-s-2} k_{i-s} f_{i-s}}
                                 {d_{i-s-1} c_{i-s-2}} 
                     \ \alpha_{i,i-s}, \ \ \ i=3,4,\ldots,n, \ \ \ s=1,2,\ldots,i-2,
\label{fra4}
\end{equation}
where $c_i$, $d_i$, $f_i$, and $g_i$ are given by the relations~(\ref{A1conventions}), computes $A_1^{-1}$ in $5 n^2 / 2 + 5 n / 2 - 6$ mult/div (since the coefficients of $\alpha_{i,i-s}$ depends only on the second subscript) and $5 n - 9$ add/sub.

In terms of $j$, the above algorithm takes the form
\begin{displaymath}
\alpha_{j-1,j} = - 1 / c_{j-1}, \ \ \ j=2,3,\ldots,n,
\end{displaymath}
\begin{displaymath}
\alpha_{jj} = - \alpha_{j,j+1} + \frac{b_{j-1} g_j}{c_{j-1} c_j}, \ \ \ 
              j=2,3,\ldots,n-1, \ \ \ \alpha_{11}=\frac{k_2}{k_1 c_1}, \ \ \ 
              \alpha_{nn} = \frac{b_{n-1}}{c_{n-1} c_n},
\end{displaymath}
\begin{displaymath}
\alpha_{j+1,j} = - \ \frac{d_{j-1} g_{j+1}}{c_{j-1} c_j c_{j+1}}, \ \ \ j=1,2,\ldots,n-1,
\end{displaymath}
\begin{displaymath}
\alpha_{j+s+1,j} = - \ \frac{g_{j+s+1} k_{j+s} f_{j+s}}
                                 {g_{j+s} c_{j+s+1}} 
                     \ \alpha_{j+s,j}, \ \ \ j=1,2,\ldots,n-2, \ \ \ s=1,2,\ldots,n\!-\!j\!-\!1.
\end{displaymath}

For the computation of $A_2^{-1}$ the algorithm~(\ref{fra1})--(\ref{fra4}) changes only in the estimation of the diagonal elements, for which we have 
\begin{displaymath}
\alpha_{ii} = - \alpha_{i,i+1} + \frac{a_{i-1} g_i}{c_{i-1} c_i}, \ \ \ i=2,3,\ldots,n-1, \ \ \
                \alpha_{11} = - \alpha_{12}, \ \ \ \alpha_{nn} = \frac{k_{n-1} b_{n-1}}{k_n c_{n-1} c_n},
\end{displaymath}
where $c_i$, $d_i$, $f_i$, and $g_i$ are given by the relations~(\ref{A2conventions}). Therefore, considering the relations~(\ref{A1conventions}) and (\ref{A2conventions}), it is clear that the number of mult/div and add/sub in computing $A_2^{-1}$ is the same with that of $A_1^{-1}$.

\section{Concluding Remarks}
The matrices $A_1$ and $A_2$ represent generalizations of known classes of test matrices. For instance, the test matrices given in \citep{mil68} (Eqs.~(2.1), (2.2)) and in \citep{her69} (Eq.~(2)) belong to the categories presented. Furthermore, by restricting the $a$'s and $b$'s to unity, $A_1$ and $A_2$ reduce to the matrices given in \citep{val77}. Also, the matrices in \citep{gka69} (pp. 41, 42, 49) are special cases of $A_1$ and $A_2$. On the other hand, concerning the recursive algorithms given in Sec.~\ref{complexity}, we have performed numerical experiments by assigning random values to the parameters of $A_1$, and with a variety of order $n$ from  $256$ to $1024$. We have found that computing $A_1^{-1}$ by the recursive algorithm~(\ref{fra1})--(\ref{fra4}) is $\sim 100$ times faster than using the LU decomposition when $n=256$ and increases gradually to $\sim 1000$ times faster when $n=1024$.


\begin{thebibliography}{}

\bibitem[Herbold~(1969)]{her69} R. J. Herbold, 
``A generalization of a class of test matrices,''
Math. Comp.,  Vol.~23, 1969, pp. 823--826.
\underline{doi:10.1090/S0025-5718-1969-0258259-0}

\bibitem[Valvi~(1977)]{val77} F. N. Valvi, 
``Explicit presentation of the inverses of some types of matrices,''
J. Inst. Maths. Applics., Vol.~19, 1977, pp. 107--117.
\underline{doi:10.1093/imamat/19.1.107} 

\bibitem[Gover \& Barnett~(1986)]{gba86} M. J. C. Gover and S. Barnett, 
``Brownian matrices: properties and extensions,''
Int. J. Systems Sci., Vol.~17, 1986, pp. 381--386.
\underline{doi:10.1080/00207728608926813} 

\bibitem[Picinbono~(1983)]{pic83} B. Picinbono, 
``Fast algorithms for Brownian matrices,''
IEEE Trans. Acoust. Speech Sig. Proc., Vol.~31, 1983, pp. 512--514.
\underline{doi:10.1109/TASSP.1983.1164078}

\bibitem[Carayannis, Kalouptsidis, \& Manolakis~(1982)]{ckm82} G. Carayannis, N. Kalouptsidis and D. G. Manolakis, 
``Fast recursive algorithms for a class of linear equations,''
IEEE Trans. Acoust. Speech Sig. Proc., Vol.~30, 1982, pp. 227--239.
\underline{doi:10.1109/TASSP.1982.1163876}

\bibitem[Milnes~(1968)]{mil68} H. W. Milnes, 
``A note concerning the properties of a certain class of test matrices,''
Math. Comp., Vol.~22, 1968, pp. 827--832.
\underline{doi:10.1090/S0025-5718-1968-0239743-1}

\bibitem[Gregory \& Karney~(1969)]{gka69} R. T. Gregory and D. L. Karney, 
``A collection of matrices for testing computational algorithms,''
Wiley--Interscience, London, 1969.











\end{thebibliography}
\end{document}